\date{Rev. 29/XII/2009 JM}
\title{The number of unit distances is almost linear for most norms}

\newcommand{\cmt}[1]{\ifhmode\newline\fi{\sf *** \ \ #1 \\}}

\author{
{\sc Ji\v{r}\'{\i} Matou\v{s}ek}
\\
   {\footnotesize Department of Applied Mathematics and}\\[-1.5mm]
   {\footnotesize Institute of Theoretical Computer Science (ITI)}\\[-1.5mm]
   {\footnotesize  Charles University, Malostransk\'{e} n\'{a}m. 25}\\[-1.5mm]
{\footnotesize  118~00~~Praha~1,
   Czech Republic}
}

\documentclass[11pt]{article}

\usepackage{a4wide}
\usepackage{latexsym}
\usepackage{amssymb}
\usepackage{epsfig}

\newtheorem{theorem}{Theorem}[section]

\newtheorem{prop}[theorem]{Proposition}

\newtheorem{lemma}[theorem]{Lemma}
\newtheorem{claim}[theorem]{Claim}

\newcommand{\heading}[1]{\vspace{1ex}\par\noindent{\bf #1}}

\makeatletter
\newcommand{\ProofEndBox}{{\ifhmode\unskip\nobreak\hfil\penalty50 \else
          \leavevmode\fi\quad\vadjust{}\nobreak\hfill$\Box$
            \finalhyphendemerits=0 \par}}
\makeatother
\newcommand{\proofend}{\ProofEndBox\smallskip}

\newcommand{\R}{{\mathbb{R}}}

\newcommand{\N}{{\mathbb{N}}}
\newcommand\eps{\varepsilon}

\newcommand\BB{\mathcal{B}}
\newcommand\MM{\mathcal{M}}

\newcommand\makevec[1]{{\bf #1}}
\def \uu {\makevec{u}}

\def \aa {\makevec{a}}
\def \bb {\makevec{b}}

\def \xx {\makevec{x}}
\def \yy {\makevec{y}}

\def \pp {\makevec{p}}
\def \zero{\makevec{0}}
\def\bG{\makevec{G}}
\def\bt{\makevec{t}}

\def\:{\colon}

\newcommand{\alterdef}[1]{\!\left\{\!\!\begin{array}{ll}
                                   #1 \end{array}  \right. }

\long\def\onefigure#1#2{
\begin{figure*}[tbp]
\begin{center}
#1
\end{center}
\caption{#2}
\end{figure*}
}

\newcommand{\labepsfig}[2]  
{\onefigure{\mbox{\epsfig{file=#1.eps}}}{\label{f:#1} #2} }

\newcommand{\labepsfigw}[3]  
{\onefigure{\mbox{\epsfig{file=#1.eps,width=#2}}}{\label{f:#1} #3} }

\newcommand{\imbal}{{\rm imb}}

\begin{document}

\maketitle

\begin{abstract} 
We prove that there exists a norm in the plane
under which no $n$-point set determines more than $O(n\log n\log\log n)$
unit distances. Actually, most norms have this property, in the sense
that their complement is a meager set in the metric space of all norms 
(with the metric given by the Hausdorff distance of the unit balls).
\end{abstract}

\section{Introduction}

What is the maximum possible number $u(n)$ of unit distances
determined by an $n$-point set in the Euclidean plane?
This tantalizing question, raised by Erd\H{o}s \cite{e-sdnp-46} in 1946,
has motivated an extensive research
(see, e.g., Brass, Moser, and Pach \cite{BMP-problems} for a survey),
but it remains wide open. 

Erd\H{o}s \cite{e-sdnp-46} proved a lower bound
$u(n)=\Omega(n^{1+c/\log\log n})$ for a constant $c>0$,
attained for the $\sqrt n\times\sqrt n$ grid, and he conjectured
that it has the right order of magnitude (and in particular,
that $u(n)=O(n^{1+\eps})$ for every fixed $\eps>0$).
However, the current best upper bound is only $O(n^{4/3})$.
It was first proved by Spencer, Szemer\'edi, and Trotter \cite{sst-udep-84},
based on the method of Szemer{\'e}di and Trotter \cite{st-cdbep-83},
and several simpler proofs are available by now
(by Clarkson et al.~\cite{cegsw-ccbac-90},
by Aronov and Sharir \cite{AronovSharir-newincid}, and the
simplest one by Sz\'ekely \cite{s-cnhep-97}).

The problem of unit distances has also been considered
for norms other than the Euclidean one. For a norm\footnote{We
recall that a (real) \emph{norm} on a real vector space $Z$ is a mapping
that assigns a nonnegative real number
$\|\xx\|$ to each $\xx\in Z$ so that
$\|\xx\|=0$ implies $\xx=0$, $\|\alpha \xx\|=|\alpha|\cdot\|\xx\|$
 for all $\alpha\in\R$,
and the triangle inequality holds: $\|\xx+\yy\|\leq \|\xx\|+\|\yy\|$.
The \emph{unit ball} of the norm $\|.\|$ is the set
$B_{\|.\|}=\{\xx\in Z: \|\xx\|\le 1\}$. The unit ball of any norm
is a closed bounded convex body $B$ that is
symmetric about $\zero$ and contains $\zero$ in the interior.
Conversely, every $B\subset Z$ with the listed properties is the unit ball
of a (uniquely determined) norm.
} 
$\|.\|$ on $\R^2$, let $u_{\|.\|}(n)$ denote the maximum possible number
of unit distances determined by $n$ points in $(\R^2,\|.\|)$.

If the boundary of the unit ball $B_{\|.\|}$
of $\|.\|$ contains a straight segment,
then it is easy to construct $n$-point sets with $\Omega(n^2)$
unit distances. On the other hand, if $\|.\|$ is 
\emph{strictly convex}, meaning that the boundary of $B_{\|.\|}$
contains no straight segment, then $u_{\|.\|}(n)=O(n^{4/3})$,
as can be shown by a straightforward generalization
of the known proofs for the Euclidean case.

Valtr \cite{valtr-norm-ud}, strengthening an earlier result of
Brass, constructed a strictly convex norm $\|.\|$ in the plane
with $u_{\|.\|}=\Omega(n^{4/3})$, thus showing that the upper
bound cannot be improved in general for strictly convex norms.

A simple construction shows that $u_{\|.\|}(n)=\Omega(n\log n)$
holds for every norm $\|.\|$ (see, e.g., \cite{BMP-problems}).
Here we will show that there exists a norm $\|.\|$ 
with $u_{\|.\|}(n)=O(n\log n\log\log n)$, almost matching the lower bound.
Actually, we show that most norms, in the sense of Baire category,
have this property.

To formulate this result, we recall the relevant notions.
Let $\BB$ be the set of all unit balls of norms in $\R^2$, i.e.,
of all closed bounded $\zero$-symmetric convex sets containing $\zero$ in
the interior. Endowed with the Hausdorff metric\footnote{We recall
that the Hausdorff distance $d_H(A,B)$ of two sets in the Euclidean
plane is defined as $\min(h(A,B),h(B,A))$, where $h(A,B)=
\sup_{a\in A}\inf_{b\in B}\|a-b\|_2$, with $\|.\|_2$ denoting the
Euclidean distance.}
$d_H$, the set $\BB$ forms a Baire space, meaning that 
each meager set\footnote{A set $S$ in
a metric (or topological) space $X$ is \emph{nowhere dense}
if every nonempty open set $U\subseteq X$ contains a nonempty
open set $V$ with $V\cap S=\emptyset$. A \emph{meager set}
is a countable union of nowhere dense sets.}
has a dense complement; see, e.g., Gruber \cite[Chapter~13]{Gruber-book}.

If $P$ is some property that a norm on $\R^2$ may or may not have,
we say that \emph{most norms have property $P$} if 
the (unit balls of the) norms not having property $P$ 
form a meager set in $\BB$. A similar terminology is 
commonly used for convex bodies. 

If most norms have
property $P_1$ and most norms have property $P_2$, then most
norms have both  $P_1$ and $P_2$ (and similarly for countably
many properties), which makes this approach a powerful tool
for proving existence results. Starting with a
paper of Klee \cite{klee-mostsmooth}, who proved
that most norms are smooth and strictly convex,
there have been many papers establishing that most norms or most convex bodies
have various properties (see \cite{Gruber-book}). 
We add the following item to this collection.

\begin{theorem}\label{t:} 
There exists a constant $C_0$
such that most norms $\|.\|$ on $\R^2$ satisfy
$$u_{\|.\|}(n) < C_0n\log n\log\log n
$$ for all $n\ge 3$
(\,$\log$ stands for logarithm in base~$2$ everywhere in this paper).
In particular, there exists a smooth and strictly convex 
norm $\|.\|$ with this property.
\end{theorem}

Since, as was mentioned above, $u_{\|.\|}(n)=\Omega(n\log n)$ for 
all norms, the bound in the theorem is tight up to the 
$O(\log\log n)$ factor. This factor comes out of a graph-theoretic
result, Proposition~\ref{p:} below, and I have no good guess 
whether it is really needed.

The proof of the theorem has two main parts. 
We begin with the first, purely graph-theoretic part
in Section~\ref{s:edgecol}. The result needed
for the rest of the proof is Proposition~\ref{p:}, 
asserting the existence of a certain subgraph
in every sufficiently dense graph with a given proper edge-coloring.
Its proof relies heavily on a similar result of 
P\v{r}\'{\i}v\v{e}tiv\'y, {\v{S}}kovro\v{n}, and the author
\cite{MatPrivSko} (but the presentation below is self-contained).

Then, in Section~\ref{s:udg} we continue with the second, geometric
part of the proof of Theorem~\ref{t:}. Very roughly speaking,
using the graph-theoretic result from the first part of the proof,
we show that if there is a set $P$ with many unit distances, under any norm,
and if $\uu_1,\ldots,\uu_k$ are all the mutually non-parallel unit vectors
defined by pairs of points of $P$, then there are ``many'' linear
dependences among the $\uu_i$. Namely, there is an integer $\ell$,
such that some $\ell+1$ vectors among the $\uu_i$
can be expressed as linear functions of some other 
$\ell$ of the $\uu_i$ (where the linear functions don't depend on the norm).
Finally, we show that most norms don't admit such linear dependences---this
is done by approximating the unit ball of the considered norm by
a convex polygon, and employing a linear-algebraic perturbation
argument to the lines bounding the polygon.


It would be interesting to prove a similar result for some narrower
class of norms. For example, one might hope to prove that
the $\ell_p$ norms admit only a near-linear
number of unit distances for most $p$ (in the Baire category sense or
even for almost all $p$ w.r.t.\ the Lebesgue measure).
For that, the idea of polygonal approximations seems unusable,
but perhaps more powerful tools from algebraic geometry might help. 

Finally, of course, it might be possible to use some pieces from
the method of this paper for attacking the Euclidean case.
However, since the number of unit distances for the Euclidean
case can be much larger than $n\log n\log \log n$, 
additional ideas are certainly needed.




\section{Connected subgraphs with few colors in edge-colored graphs}
\label{s:edgecol}


Let $G=(V,E)$ be a (simple, undirected) graph. 
An \emph{edge coloring} of $G$ is a mapping $c\:E\to\N=\{1,2,3,\ldots\}$.
The edge coloring $c$ is called \emph{proper} if $c(e)\ne c(e')$
whenever the edges $e$ and $e'$ share a vertex.

Let $G$ be a graph with a given edge coloring.
For a subset $W\subseteq V$ of vertices we let $G[W]$ stand for
the subgraph of $G$ induced by $W$, with the edge 
coloring inherited from that of $G$. Further, if $I\subseteq \N$
is a set of colors, we write $G[I,W]$ for the
subgraph induced by $W$ on the edges with colors in $I$, that is,
$$
G[I,W]=\Bigl(W,\{\{u,v\}\in E: u,v\in W,c(e)\in I\}\Bigr)
$$
(the coloring is not explicitly mentioned in the notation).

\begin{prop}\label{p:} 
Let $q>1$ be a real parameter.
Let $G=(V,E)$ be a graph on $n\ge 4$ vertices,
with at least $Cq n\log n\log\log n$ edges (where $C$
is a suitable absolute constant), and with a given proper edge
coloring. Then there exist
a nonempty subset $W\subseteq V$ of vertices, $|W|\ge 2$, and a subset
$I\subset\N$ of colors such that the subgraph $G[I,W]$
is connected and the edges of $G[W]$ have at least
$q|I|$ distinct colors.
\end{prop}

As was mentioned in the introduction, this proposition is
similar to a result from \cite{MatPrivSko}, and the proof
is also quite similar to the one in \cite{MatPrivSko}.
I still consider it worth presenting in full, since
describing the required modifications would be clumsy,
and moreover, the proof below is significantly 
simpler than that in \cite{MatPrivSko}, mainly because
the required result is weaker (in Proposition~\ref{p:}
we obtain a single connected subgraph, while in \cite{MatPrivSko}
several color-disjoint connected subgraphs on the same vertex
set were needed).
\medskip

At the beginning of the proof, we use a well-known observation
stating that every graph of average degree $\delta$
has a subgraph whose minimum degree is at least $\delta/2$
(this follows by repeatedly deleting vertices of degree
below $\delta/2$ and checking that the average degree
can't decrease). So we may assume that $G$ has minimum
degree at least $Cq n\log n\log\log n$.

Let $W\subseteq V$ be a subset of vertices of $G$ (so far arbitrary).
An \emph{edge cut} in $G[W]$ is a partition $(A,B)$ of $W$ into two
nonempty subsets. We define the \emph{maximum degree} $\Delta(A,B)$ 
of such an edge cut as the maximum number of neighbors
of a vertex from $A$ in $B$ or of a vertex from $B$ in $A$;
formally,
$$
\Delta(A,B):=\max\left\{\max_{a\in A}|\{\{a,b\}\in E:b\in B\}|,
                         \max_{b\in B}|\{\{b,a\}\in E:a\in A\}|\right\}.
$$

The proof of Proposition~\ref{p:}
proceeds in two stages. In the first stage, we forget about
the edge colors; we select 
the set $W$ so that every edge cut in $G[W]$ has a sufficiently
large maximum degree.
In order to get the (almost tight) quantitative result in the proposition,
we need to quantify the ``sufficiently large maximum degree''
of a cut depending on the \emph{imbalance} of the cut, which
is defined by
$$
\imbal(A,B) := \frac{|A|+|B|}{\min(|A|,|B|)}.
$$

\begin{lemma}\label{l:coloredcuts} Let $r\ge 1$ be a parameter
(which we will later set to $Cq\log\log n$
in the application of the lemma), and let $G=(V,E)$ be a graph
on $n\ge 2$ vertices of minimum degree at least $
r\log n$. Then there exists $W\subseteq V$,
$|W|\ge 2$,
such that every edge cut $(A,B)$ in $G[W]$ satisfies
$$
\Delta(A,B)\ge r \log \imbal(A,B).
$$
\end{lemma}

\heading{Proof. } 
The proof proceeds by a recursive partitioning: 
As long as we can find an edge cut $(A,B)$ of 
small maximum degree in the current
graph, we discard the \emph{larger} of the sets $A,B$. 

More formally, we set $V_1:=V$.
If $G[V_j]$ has already been 
constructed and if there is an edge cut
$(A_j,B_j)$ in $G[V_j]$ with $\Delta(A_j,B_j)<r\log\imbal(A_j,B_j)$,
we let $V_{j+1}$ be the smaller of the sets $A_j$ and $B_j$ (ties broken
arbitrarily) and iterate.
If there is no such edge cut, we set $W:=V_j$, $t:=j$, and finish.

It remains to
show that the resulting $W$ is nontrivial, i.e., $|W|\ge 2$.
This is clear for $t=1$ (no partition step was made),
so we assume $t\ge 2$. We show that $G[W]=G[V_t]$
has minimum degree at least~$1$, and thus $W$ can't consist
of a single vertex.

Initially, in $G$, each vertex has degree at least $r\log n$,
and by passing from $V_j$ to $V_{j+1}$, each vertex of $V_{j+1}$
loses at most $\Delta(A_j,B_j)<r\log\imbal(A_j,B_j)$ neighbors.
Thus, the minimum degree in $G[V_t]$ is strictly larger than
\begin{eqnarray*}
  r\log n - r\sum_{j=1}^{t-1} \log\imbal(A_j,B_j)
&=&r\log n - r\sum_{j=1}^{t-1}\log\frac{|V_j|}{|V_{j+1}|}\\
&=&r\log  n - r\Bigl(\log |V_1| - \log |V_{t}|\Bigr)\ge 0.\\
\end{eqnarray*}
The lemma is proved.
\proofend

\medskip

Now we continue with the \emph{second stage of 
the proof of Proposition~\ref{p:}}. Only here we start considering the
edge colors.

According to Lemma~\ref{l:coloredcuts},
we now assume that $W\subseteq V$, $|W|\ge 2$, is such that every edge cut
$(A,B)$ in $G[W]$ has maximum degree at least $r\log \imbal(A,B)$,
with $r=Cq\log\log n$. 
Consequently, the edges of every edge cut
$(A,B)$ have at least $r\log \imbal(A,B)$ distinct colors
(since the edge coloring is proper), and this
is the only property of $G[W]$ we will use.

Let $k$ denote the number of colors occurring on the edges of $G[W]$.
We note that $k\ge r$ (this follows by
using the condition above for an arbitrary cut).
It remains to show that $G[W]$ has a connected subgraph that uses at
most $k/q$ colors.

We select the colors greedily one by one, as follows.
We set $I_0:=\emptyset$, and for $j=0,1,2,\ldots,$
we do the following: If $G[I_{j},W]$ is connected, we
set $I:=I_{j}$ and finish. Otherwise, we let $i_j$
be a color $i$ minimizing the number of connected
components of $G[I_j\cup\{i\},W]$.
Then we set $I_{j+1}:=I_j\cup\{i_j\}$,
and we continue with the next step.
We need to show that we obtain a connected graph before exhausting
more than $k/q$ colors. 

Let $m_j$ be the number of connected components
of $G[I_j,W]$.  We want an upper bound on the smallest $j$ with $m_j=1$.
First we observe that $m_{j+1}\le m_j-1$ for all $j$,
since every edge cut contains at least one color. 
In the sequel, we will actually estimate
the smallest $j$ such that $m_j\le 3$.
Then at most two more steps suffice
to get down to $m_j=1$. 

We now want to bound $m_{j+1}$ in terms of $m_j$.
Essentially, we will see that adding a random color to $I_j$
is likely to connect up many components.
 
Let $K_1,\ldots,K_{m_j}$ be the vertex sets of the
connected components of $G[I_j,W]$. 
The average number of vertices in a component is $m/m_j$;
we call a component \emph{small} if it has at most $2m/m_j$
vertices. By Markov's inequality, there are at least $m_j/2$ small components.

Let $i$ be one of the colors occurring on the edges of $G[W]$ but
not belonging to $I_j$ (so there are $k-j$ possible choices for $i$).
We say that a component $K_s$ \emph{gets connected} by $i$
if there is an edge of color $i$ connecting a vertex of $K_s$
to a vertex outside $K_s$.

By the condition on the edge cuts of $G[W]$, if $K_s$ is a small
component, then the number of colors $i$ by which $K_s$ gets connected
is at least
$$
r\log \imbal\left(K_s,W\setminus K_s\right)\ge
r\log \frac m{2m/m_j}=
r\log (m_j/2).
$$
Thus, the expected number of small components that get connected
by a random color is at least
$$
\frac{m_j}{2}\cdot \frac{r\log (m_j/2)}{k-j}\ge
\frac{m_j}{2}\cdot \frac{r\log (m_j/2)}{k}.
$$
So at least this many components get connected by the color $i_{j+1}$.

It is easy to check that the number of components
always decreases at least by half of the  number of components that get
connected (an extremal case being components merged in pairs).
Thus, we have
$$
m_{j+1}\le m_j - \frac{m_j}{4}\cdot \frac{r\log (m_j/2)}{k}
\le m_j\left(1-\frac{r\log (m_j/2)}{4k}\right)
\le m_je^{-r\log (m_j/2)/4k}
$$
(we used $1-x\le e^{-x}$ in the last step). Assuming, as we may,
that $m_j\ge 4$, we have $\log (m_j/2)\ge \frac 12\log  m_j\ge
\frac12\ln m_j$, and
so 
$$
\ln m_{j+1}\le \ln m_j-r(\ln m_j)/8k=(1-r/8k)\ln m_j\le e^{-r/8k}\ln m_j.
$$
Since $m_1\le n$, we can see that $m_j$ drops below $4$ in
at most $O((k/r)\log\log n)=O(k/Cq)$ steps. We need at most two
extra colors to get all the way to $m_j=1$, so altogether the number
of colors needed to build a connected
graph is $O(k/Cq+2)=O(k/Cq)$ (since $k\ge r$, and thus $k/Cq\ge
\log\log n\ge 1$).  The implicit constant
in the $O(.)$ notation is independent of $C$, and thus we can
set $C$ so large that the number of colors is at most $k/q$. 
Proposition~\ref{p:} is proved.
\proofend

\section{Unit-distance graphs}\label{s:udg}

Let $\|.\|$ be a norm in the plane, and let 
$P=(\pp_1,\ldots,\pp_n)$ be a sequence of $n$ distinct points in the plane.
With these objects we associate a finite combinatorial object,
which we will call the \emph{decorated unit-distance graph}.

First, we define the \emph{unit-distance graph}  $G=G(\|.\|,P)$
as the (undirected) graph $(V,E)$ with vertex set $V:=[n]$ 
(where we use the notation $[n]=\{1,2,\ldots,n\}$)
and with edges corresponding to the pairs of points of unit distance; that is,
$E=\{\{a,b\}:\|\pp_b-\pp_a\|=1\}$.

To every edge $e=\{a,b\}\in E$ we assign a vector
$\uu(e)$, in such a way that $\uu(e)=\pm(\pp_b-\pp_a)$, and the sign
is chosen using some globally consistent rule, so that
parallel edges get the same $\uu(e)$. For example,
we may require that $\uu(e)$ lie in the closed upper halfplane
minus the negative $x$-axis.

Let $U:=\{\uu(e):e\in E\}$ be the \emph{unit direction set} of $P$, 
and we let $\uu_1,\uu_2,\ldots,\uu_k$
be an enumeration of all distinct elements of $U$, 
say in the lexicographic order. We call $\uu_1,\ldots,\uu_k$ the
\emph{unit directions} of $P$ (under $\|.\|$).
Then we define a coloring $c\:E\to[k]$
of the edges of the unit-distance graph, setting $c(e)=i$ if $\uu(e)=\uu_i$.
(We note that $c$ need \emph{not} be a proper
edge coloring, since there can be two edges with the same direction
incident to a single vertex.)

Finally, we record the geometric orientation of each edge.
Namely, we define a mapping $\sigma\:E\to\{-1,+1\}$:
For an edge $\{a,b\}\in E$ with $a<b$ we set
$$
\sigma(\{a,b\})=\alterdef{+1&\mbox{ if $\uu(e)=\pp_b-\pp_a$},\\
                          -1&\mbox{ if $\uu(e)=\pp_a-\pp_b$.}}
$$
The \emph{decorated unit-distance graph} of $P$ under $\|.\|$
is  defined as the triple $\bG=\bG(\|.\|,P):=(G,c,\sigma)$.

Now we define an \emph{abstract decorated unit-distance graph} as
expected, i.e., as a triptuple $\bG=(G,c,\sigma)$, 
where $G$ is a graph with vertex
set $[n]$ for some $n$, $c$ is a mapping $E(G)\to[k]$
for some $k$, and $\sigma$ is a mapping $E\to\{-1,+1\}$.
 We say that a sequence $P$ of distinct
points in $\R^2$ is a \emph{realization}
of an abstract decorated unit-distance graph $\bG$ under $\|.\|$ if 
$\bG$ is equal to the decorated unit-distance graph of $P$ under $\|.\|$.
(We require equality to keep the definitions simple; we could
as well introduce a suitable notion of isomorphism, but there is no need.)

Here is the main result of this section. Roughly speaking, it tells
us that if $\bG$ is a sufficiently dense 
abstract decorated unit-distance graph, then for every realization,
the unit directions satisfy certain fixed linear dependences---some
$\ell+1$ of the unit directions can be expressed using
some other $\ell$ of the unit directions.

\begin{lemma}\label{l:lindep}
The following holds for a sufficiently large constant $C_0$.
Let $\bG$ be an abstract decorated unit-distance graph with $n\ge 4$
vertices, at least $f(n):=C_0n\log n\log\log n$ edges, and $k$ colors. 
Then there exists an integer $\ell\ge 1$, a sequence 
$(i(1),i(2),\ldots,i(2\ell+1))$ of distinct indices in $[k]$, and 
linear maps $L_1,L_2,\ldots,L_{\ell+1}\:(\R^2)^\ell\to\R^2$ such
that for every realization $P$ of $\bG$ (under any norm),
we have 
$$
\uu_{i(\ell+j)}=L_j(\uu_{i(1)},\uu_{i(2)},\ldots,\uu_{i(\ell)}),
\ \ \ j=1,2,\ldots,\ell+1,
$$
where $\uu_1,\ldots,\uu_k$ are the unit directions of~$P$.
\end{lemma}

\heading{Proof. } Let $\bG=(G,c,\sigma)$. 
In order to apply Proposition~\ref{p:}, we may need to prune
the graph so that $c$ becomes a proper edge coloring.
If $\bG$ has any realization at all, then, for geometric
reasons, no color occurs on more
than two edges incident to each vertex. Hence, for each $i$, the subgraph
made of edges of color $i$ consists of paths and cycles, and
so by deleting at most $\frac 23$ of the edges, we can turn
this subgraph into a matching, and hence obtain a subgraph
$\tilde G$ of $G$ with at least $\frac1{3}f(n)$ edges
for which $c$ is a proper edge coloring. (By using more geometry,
it is easily seen that it even suffices to delete only at most $\frac12$
of the edges, rather than $\frac23$.)

Now we are ready to apply Proposition~\ref{p:} on the graph $\tilde G$
with the proper edge coloring $c$, and with $q=2.001$, say.
This yields a subset $W\subseteq V(\tilde G)$ and a subset $I\subset[k]$
of colors, such that the subgraph $\tilde G[I,W]$ is connected,
and $\tilde G[W]$ uses at least $2|I|+1$ colors. Let $J$ be a set
of $|I|+1$ colors used on the edges of $W$ but not belonging to $I$.

Now we can define the objects whose existence is claimed in the lemma.
We set $\ell:=|I|$, let $(i(1),\ldots,i(\ell))$ be an enumeration of $I$,
and let $i(\ell+1),\ldots,i(2\ell+1)$ be an enumeration of~$J$.

Let us consider some color $j\in J$, and let $\{a,b\}$ be an edge
of color $j$ in $\tilde G[W]$. Then there is a path $\pi$ from $a$ to $b$
in $\tilde G[W]$
whose edges have only colors in $I$,
and for every realization $P$ of $\bG$, $\uu_j$ is a signed sum
of the unit directions along this path. An example is given
in Fig.~\ref{f:sigsum}: If $a=1$, $b=6$, $j=1$, the edge $\{1,6\}$
has sign $-1$, the path
$\pi$ goes through the vertices $2,3,4,5$ in this order,
and its edges have colors $2,2,3,4,2$ and signs $+1,+1,-1,+1,+1$, then
$\uu_1=-3\uu_2+\uu_3-\uu_4$. 
This yields the desired linear maps $L_1,\ldots,L_{\ell+1}$,
and the lemma is proved. \proofend

\labepsfig{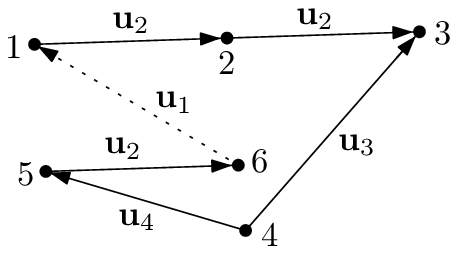}{Expressing $\uu_j$ in terms of the $\uu_i$, $i\in I$.}

\section{Proof of Theorem~\ref{t:}}

Let us call a norm $\|.\|$ on $\R^2$ \emph{bad} if
$u_{\|.\|}(n)\ge f(n)=C_0n\log n \log\log n$ for some $n\ge 3$,
and let $\MM\subseteq\BB$ be the set of all bad norms.
We want to show that $\MM$ is meager, and thus we want
to cover it by countably many nowhere dense sets.

In our proof, the nowhere dense sets $\MM_{\bG,\eta}$ are indexed by two
parameters: $\bG$, which runs through all abstract decorated
unit-distance graphs with $n$ vertices and at least
$f(n)$ edges, $n=3,4,\ldots$,
and $\eta$, which runs through all positive
numbers of the form $\frac 1m$, $m$ an integer.

To define $\MM_{\bG,\eta}$, we first define that a realization
$P$ of $\bG$ under a norm $\|.\|$ is \emph{$\eta$-separated}
if for every two unit direction vectors $\uu_i,\uu_j$ of this
realization, the lines spanned by $\uu_i$ and $\uu_j$ 
have angle at least $\eta$.

Now $\MM_{\bG,\eta}$ consists of all norms $\|.\|$ under which
$\bG$ has an $\eta$-separated realization. 

It is easily checked that
the $\MM_{\bG,\eta}$ cover all of $\MM$. Indeed, for every bad norm
$\|.\|$ we can choose $n$ and an $n$-point sequence $P$
with at least $f(n)$ unit distances. We define $\bG$
as the decorated unit-distance graph of $P$ under $\|.\|$.
It remains to observe that, trivially, every realization
of $\bG$ under some norm is $\eta$-separated for some $\eta>0$.
Thus $\|.\|\in\MM_{\bG,\eta}$.

The main part of the proof consists of showing that each $\MM_{\bG,\eta}$
is nowhere dense. Explicitly, this is expressed in the following lemma;
once we prove it, we will be done with Theorem~\ref{t:}
(the smoothness and strict convexity asserted in the theorem
follows from Klee's result \cite{klee-mostsmooth}
mentioned in the introduction, namely,
that most norms are smooth and strictly convex).

\begin{lemma}\label{l:nowd}
 Let $\bG$ be an abstract decorated unit-distance graph
with $n$ vertices and at least $f(n)$ edges, let $B_0\in \BB$ be 
the unit ball of some norm, and let $\eta,\eps>0$. Then there
exist $B\in \BB$ with $d_H(B,B_0)<\eps$
(where $d_H$ denotes the Hausdorff distance) 
and $\delta>0$ such that no $B'\in \BB$ with $d_H(B',B)<\delta$
belongs to $\MM_{\bG,\eta}$.
\end{lemma}

\heading{Proof. } First we approximate $B_0$ by a $\zero$-symmetric convex
polygon $B_1$  within Hausdorff distance at most $\frac\eps2$ from $B_0$.
We make sure that all sides of $B_1$ are sufficiently short,
so short that two lines through $\zero$ with angle at least $\eta$
never meet the same side 
of $B_1$.
(If $B_0$ has straight segments in the boundary, we need to
to ``bulge'' $B_1$ slightly; see Fig.~\ref{f:polapprox}.)

\labepsfig{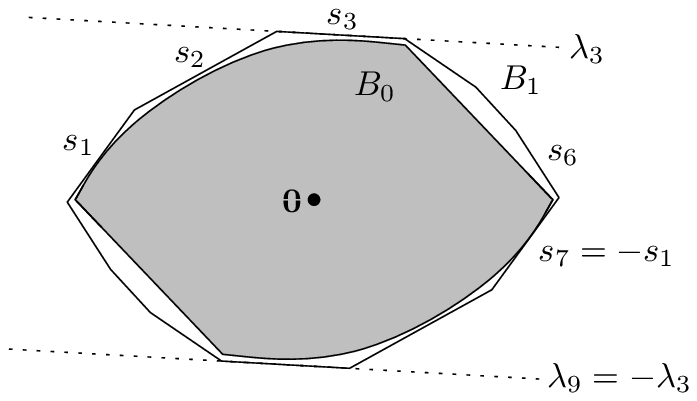}{Approximating the unit ball $B_0$ by a convex
polygon. }

Let $s_1,s_2,\ldots,s_{2m}$ be the sides of $B_1$ listed in clockwise
order, say, so that $s_{i}$ and $s_{m+i}$ are opposite
(i.e., $s_{m+i}=-s_i$).
Let $\lambda_i$ be the line spanned by $s_i$, and for a real parameter
$t$, let $\lambda_i(t)$ be the line obtained by a parallel translation
of $\lambda_i$ by distance $t$, where $t>0$ means translation
 away from the origin and $t<0$ translation towards the origin. 
We have $\lambda_{m+i}(t)=-\lambda_i(t)$.

Let us consider an $m$-tuple
$\bt=(t_1,\ldots,t_m)\in T_0:=[-\delta_0,\delta_0]^m$. 
For $\delta_0>0$ sufficiently
small, the lines $\lambda_1(t_1),\ldots,\lambda_m(t_m),\lambda_{m+1}(t_1),
\ldots,\lambda_{2m}(t_m)$ bound a symmetric convex polygon with $2m$ sides,
which we denote by $B_1(\bt)$. Moreover, for $\delta_0$ sufficiently small,
$d_H(B_1(\bt),B_0)<\eps$,
and the sides of $B_1(\bt)$ are still short in the same sense as those of~$B_1$.

Now we digress from geometry for a moment and we apply Lemma~\ref{l:lindep}
to the abstract decorated unit-distance graph $\bG$. This yields an integer
$\ell$, indices $i(1),\ldots,i(2\ell+1)$, and linear maps
$L_1,\ldots,L_{\ell+1}$ as in the lemma. In order to make the notation
slightly simpler, let us pretend that $i(j)=j$ for all $j=1,\ldots,2\ell+1$.
Thus, for every realization of $\bG$, the unit directions 
$\uu_1,\ldots,\uu_{2\ell+1}$ satisfy the linear relations
$\uu_{\ell+i}=L_i(\uu_1,\ldots,\uu_\ell)$, $i=1,2,\ldots,\ell+1$.

Next, let us consider a particular realization of $\bG$ under the norm
induced by $B_1(\bt)$ for  some $\bt\in T_0$. Each of the unit directions
$\uu_1,\ldots,\uu_{2\ell+1}$ lies on the boundary of $B_1(\bt)$,
and thus on some line $\lambda_\alpha(t_\alpha)$. (Here we abuse the notation
slightly, since the range of $\alpha$ is $[2m]$, while
$\bt$ is indexed only by $[m]$, in order to preserve the symmetry
of the polygon. So we make the convention that $t_{m+i}$ is the same
as $t_i$.)

Let $\alpha(i)\in [2m]$ 
be the index
such that $\uu_i$ lies on $\lambda_{\alpha(i)}(t_{\alpha(i)})$ 
(if $\uu_i$ is a vertex of
the polygon and thus lies on two of the lines, we pick one arbitrarily).
Since the sides of $B_1(\bt)$ are short,
we have $\alpha(i)\ne\alpha(i')$ whenever $i\ne i'$,
and also $\alpha(i)+m)\ne\alpha(i')$ (where $\alpha(i)+m$ is to
be understood modulo $2m$).

Let us call a mapping $\alpha\:[2\ell+1]\to[2m]$ an \emph{admissible
assignment of lines} if it satisfies the condition in the previous
sentence. Let us define a \emph{box} $T\subseteq T_0$
as a product of closed intervals with a nonempty interior;
each box can be written as an $m$-dimensional ``interval''
$[\bt_{\rm min},\bt_{\rm max}]$.
Our next goal is establishing the following claim.

\begin{claim}\label{c:laim} 
There exists a box $\tilde T\subseteq T_0$ such that
and for every admissible assignment
of lines $\alpha$ and for every $\bt\in \tilde T$
there are no vectors $\uu_1,\ldots,\uu_{2\ell+1}\in\R^2$
such that each $\uu_i$ lies on the appropriate
line, i.e., $\uu_i\in\lambda_{\alpha(i)}(t_{\alpha(i)})$,
and the $\uu_i$ satisfy the linear relations
$\uu_{\ell+i}=L_i(\uu_1,\ldots,\uu_\ell)$, $i=1,2,\ldots,\ell+1$.
\end{claim}

\heading{Proof of the claim. } We will kill all admissible
assignments $\alpha$ one by one inductively,
progressively shrinking the current box. The following
statement allows us to make an inductive step:
\emph{Let $T\subseteq T_0$ be a box,
and let $\alpha$ be an admissible assignment of lines. 
Then there exists a box $T'\subseteq T$ such that for every $\bt\in T'$
there are no vectors $\uu_1,\ldots,\uu_{2\ell+1}\in\R^2$
with $\uu_i\in\lambda_{\alpha(i)}(t_{\alpha(i)})$ for all $i$
and with $\uu_{\ell+i}=L_i(\uu_1,\ldots,\uu_\ell)$, $i=1,2,\ldots,\ell+1$.}

To prove this, let us consider a vector $\xx\in\R^{2\ell}$,
which we think of as a concatenation of $\uu_1,\ldots,\uu_{\ell}$,
and let us think of its components $x_i$ as unknowns.

For each $i=1,2,\ldots,\ell$, the condition 
$\uu_i\in \lambda_{\alpha(i)}(t_{\alpha(i)})$ translates to a single linear
equation for $\xx$, of the form $\aa_i^T\xx=b_i$, where
the coefficient vector $\aa_i$ on the left-hand side doesn't depend
on $\bt$, while $b_i=b_i(t_{\alpha(i)})$ is a \emph{nonconstant} linear
function of $t_{\alpha(i)}$.

Similarly, for $i=1,2,\ldots,\ell+1$, the condition
$\uu_{\ell+i}\in \lambda_{\alpha(\ell+i)}(t_{\alpha(\ell+i)})$
together with $\uu_{\ell+i}=L_i(\uu_1,\ldots,\uu_\ell)$ translate
to a similar linear equation $\aa_{\ell+i}^T\xx=b_{\ell+i}$,
again with $\aa_{\ell+i}$ independent of $\bt$ and with
$b_{\ell+i}=b_{\ell+i}(t_{\alpha(\ell+i)})$  a nonconstant linear
function of $t_{\alpha(\ell+i)}$. 

Since the $\alpha(i)$ are all
distinct, altogether we get that if the appropriate $\uu_1,\ldots,\uu_{2\ell+1}$
exist, then $\xx$ satisfies the system $A\xx=\bb$ of $2\ell+1$ linear
equations with $2\ell$ unknowns, where $A$ is a fixed matrix
and the right-hand side $\bb=\bb(\bt)$ is a \emph{surjective}
linear function $\R^{m}\to\R^{2\ell+1}$.

Since we have more equations than unknowns, the system $A\xx=\bb$
has a solution only for $\bb$ contained in  a proper linear
subspace of $\R^{2\ell+1}$. Hence,
by the surjectivity of $\bb(\bt)$, the set of all $\bt\in\R^{m}$ 
for which  $A\xx=\bb(\bt)$ is unsolvable is a dense open
subset of $\R^{m}$. From this the existence of the desired box $T'$
follows, and the Claim~\ref{c:laim} is proved.
\proofend

\heading{Finishing the proof of Lemma~\ref{l:nowd}. }
Let us consider the box $\tilde T=[\bt_{\rm min},\bt_{\rm max}]$ 
as in Claim~\ref{c:laim}. We set $\bt_{\rm mid}:= 
(\bt_{\rm min}+\bt_{\rm max})/2$, and we consider the polygons
 $B_{\rm in}:=B_1(\bt_{\rm min})$, $B_{\rm out}:=B_1(\bt_{\rm max})$,
and $B:=B_1(\bt_{\rm mid})$; see Fig.~\ref{f:polygs}.
We claim that $B$ is as in the lemma, i.e., no $B'\in\BB$ sufficiently
close to $B$ belongs to $\MM_{\bG,\eta}$.
\labepsfig{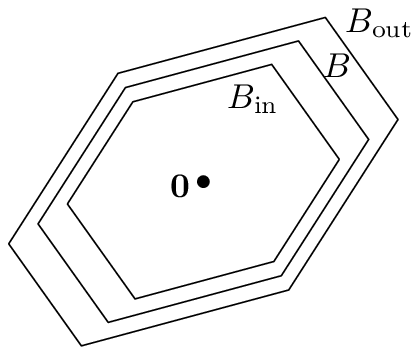}{The polygons $B_{\rm in}$, $B_{\rm out}$, and $B$.}

To see this, we note that every $B'$ sufficiently close
to $B$ satisfies $B_{\rm in}\subseteq B'\subseteq B_{\rm out}$.
For contradiction, we assume that there is an $\eta$-separated
realization of $\bG$ under $B'$. Then the unit directions
$\uu_1,\ldots,\uu_{2\ell+1}$ lie on the boundary of $B'$.

The region $B_{\rm out}\setminus B_{\rm in}$ is naturally divided
into $2m$ trapezoids $R_1,\ldots,R_{2m}$ belonging to the sides,
as in Fig.~\ref{f:polygs1}. Each of $\uu_i$, $i=1,2,\ldots,2\ell+1$,
lies in one of these trapezoids, let us call it $R_{\alpha(i)}$
(border disputes resolved arbitrarily).
Since the considered realization is $\eta$-separated,
no two of the $\uu_i$ share the same trapezoid, and also no two
of these trapezoids are opposite to one another. 
So $\alpha$ defines an admissible assignment of sides.

\labepsfig{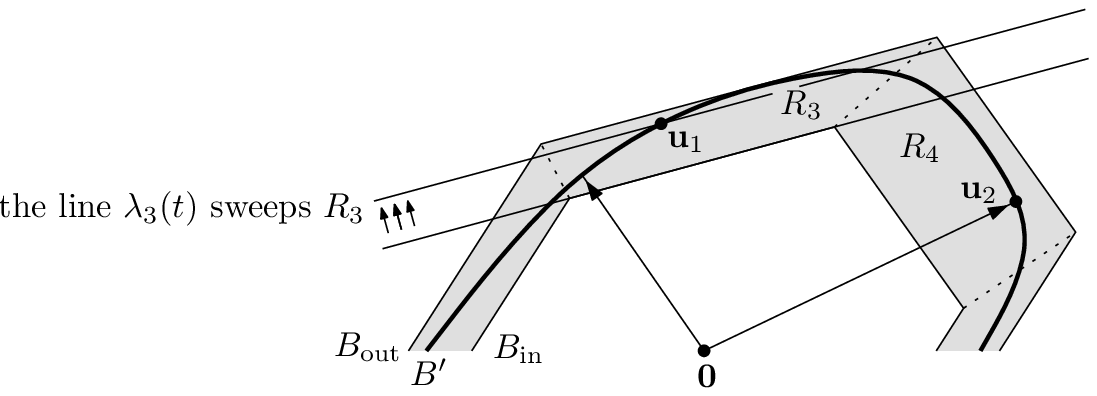}{Dividing the region $B_{\rm out}\setminus B_{\rm in}$
into trapezoids.}

Let us consider the trapezoid $R_{\alpha(i)}$.
As the line $\lambda_{\alpha(i)}(t)$ moves from the inner position
(with $t=(\bt_{\rm min})_{\alpha(i)}$) to the outer position
(with $t=(\bt_{\rm max})_{\alpha(i)}$), it sweeps the whole of $R_{\alpha(i)}$,
and hence for some $t$ it
contains $\uu_i$; let us denote this value of $t$ by $\bar t_{\alpha(i)}$.

This defines $2\ell+1$ of the components of a vector $\bar\bt\in\R^m$.
Let us set the remaining components to the corresponding components
of $\bt_{\rm mid}$, say. 
Then $\bar\bt$
lies in the box $\tilde T$, and hence, by Claim~\ref{c:laim},
 $\uu_1,\ldots,\uu_{2\ell+1}$
cannot lie on the corresponding lines. The resulting contradiction
proves the lemma, and this also finishes the proof of Theorem~\ref{t:}. 
\proofend


\newcommand{\etalchar}[1]{$^{#1}$}

\end{document}